\documentclass[a4paper, 12pt]{article}

\usepackage[utf8]{inputenc}
\usepackage[font=small,labelfont=bf]{caption}
\usepackage[skip=0.4\baselineskip]{caption}
\usepackage{rotating}
\usepackage{longtable}
\usepackage{amsmath,amsfonts,amssymb}
\usepackage{mathtools}

\usepackage{multirow}
\usepackage{fullpage}
\usepackage{times}
\usepackage{adjustbox}
\usepackage{pgfplots}
\usepackage{tablefootnote}
\usepackage{fancyhdr,graphicx,amsmath,amssymb}
\include{pythonlisting}
\usepackage[english]{babel}
\usepackage[T1]{fontenc}
\usepackage{placeins}
\usepackage{listings}

\usepackage{setspace}
\usepackage{amssymb}
\usepackage{bbm}
\usepackage{algorithm}
\usepackage{algpseudocode}
\usepackage{accents}

\usepackage{makecell}
\usepackage{amsthm}

\theoremstyle{plain}
\newtheorem{theorem}{Theorem}[section]

\theoremstyle{definition}
 
\theoremstyle{remark}

\theoremstyle{property}
\newtheorem{property}[theorem]{Property}

\algnewcommand{\algorithmicor}{\textbf{ or }}
\algnewcommand{\OR}{\algorithmicor}

\usepackage[a4paper,margin=1in]{geometry}
\usepackage{eurosym}
\usepackage{amsmath}
\usepackage{natbib}
\usepackage{graphicx}
\usepackage[colorinlistoftodos]{todonotes}
\usepackage[colorlinks=true, allcolors=blue]{hyperref}
\usepackage{authblk}
\usepackage{cleveref}
\usepackage[super]{nth}
\usepackage{enumitem, float, subcaption}
\usepackage[compact]{titlesec}
\usepackage{booktabs,caption}
\usepackage[flushleft]{threeparttable}

\newcommand{\matr}[1]{\mathbf{#1}}

\newcommand{\RNum}[1]{\uppercase\expandafter{\romannumeral #1\relax}}

\makeatletter
\newenvironment{smallalign*}
{\par
\scriptsize
\setlength{\abovedisplayskip}{0pt} 
\csname align*\endcsname}
{\csname endalign*\endcsname}
\makeatother

\usepackage{changes}
\usepackage{lipsum}
\definechangesauthor[name={Rick Willemsen}, color=orange]{rick}

\begin{document}
\hypersetup{
    colorlinks = true, 
    citecolor = black,
    linkcolor = black,
    urlcolor = black
}

\title{Generating Random Vectors satisfying Linear and Nonlinear Constraints}
\author{Rick S. H. Willemsen\thanks{Corresponding author: \href{mailto:willemsen@ese.eur.nl}{willemsen@ese.eur.nl}}}
\author{Wilco van den Heuvel}
\author{Michel van de Velden}
\affil{Econometric Institute, Erasmus University, 3062 PA Rotterdam, The Netherlands}

\maketitle

\renewenvironment{abstract}
{\begin{quote}
\noindent \rule{\linewidth}{.5pt}}
{\smallskip\noindent \rule{\linewidth}{.5pt}
\end{quote}
}

\begin{abstract} \\
We consider the problem of generating $n$-dimensional vectors with a fixed sum, with the goal of generating a uniform distribution of vectors over a valid region. This means that each possible vector has an equal probability of being generated. The Dirichlet-Rescale (DRS) algorithm, introduced by \cite{griffin2020}, aims to generate a uniform distribution of vectors with fixed sum that satisfies lower and upper bounds on the individual entries. However, we demonstrate that the uniform distribution property of the DRS algorithm does not hold in general. Using an analytical procedure and a statistical test, we show that the vectors generated by the DRS algorithm do not appear to be drawn from a uniform distribution. To resolve this issue, we propose the Dirichlet-Rescale-Constraints (DRSC) algorithm, which handles more general constraints, including both linear and nonlinear constraints, while ensuring that the vectors are drawn from a uniform distribution. In our computational experiments we demonstrate the effectiveness of the DRSC algorithm.
\\
\end{abstract}
\vfill

\newpage
\section{Introduction}
We study the generation of vectors with a fixed sum that satisfy linear and nonlinear constraints. It is common to impose that the vectors should be drawn from a uniform distribution, meaning that all vectors satisfying the constraints should be generated with an equal probability. Methods that generate vectors with a fixed sum are commonly applied to generate so-called utilisation vectors, which are used in the experimental evaluation of schedulability tests, see for instance \cite{griffin2020}.

Let $\matr{x}$ be a random $n$-dimensional vector, where we assume without loss of generality that the entries of $\matr{x}$ sum to 1. {Efficient procedures for sampling vectors with fixed sum from a uniform distribution are described in \cite{bini2005} and \cite{davis2011}.} \cite{stafford2006} proposed a method for generating vectors $\matr{x}$ satisfying global lower and upper bounds on the entries of vector $\matr{x}$, denoted by $l$ and $b$, respectively. Specifically, this ensures that $l\matr{1} \leq \matr{x} \leq b\matr{1}$, where $\matr{1}$ is a vector of ones of appropriate length. This method has been published in \cite{emberson2010}. \cite{griffin2020} proposed the Dirichlet-Rescale (DRS) algorithm that enforces lower and upper bounds for individual entries, denoted by vectors $\matr{l}$ and $\matr{u}$, respectively. The aim of the DRS algorithm is to generate vectors $\matr{x}$ from a uniform distribution satisfying $\matr{l} \leq \matr{x} \leq \matr{u}$. \cite{griffin2020} conduct a statistical test to compare the DRS algorithm with the method proposed by \cite{davis2011}, which has been proven to produce a uniform distribution. Based on a large number of generated vectors, \cite{griffin2020} found no evidence to reject the null hypothesis, which suggested that the DRS algorithm generates a uniform distribution of vectors.

However, we demonstrate that the DRS algorithm does not necessarily generate vectors from a uniform distribution. This may occur when multiple bounds are violated at the same time. To prove this, we make use of both an analytical procedure and a statistical test on specifically constructed examples. First, we analyse a simplified version of the DRS algorithm with the goal of generating $3$-dimensional vectors. We provide an analytical procedure to demonstrate that the vectors obtained by the simplified DRS are not generated from a uniform distribution. Second, we consider the original DRS algorithm as implemented by \cite{griffin2020}. We presume that $3$-dimensional vectors generated by the original DRS algorithm are from a uniform distribution, as an unintentional consequence of an additional transformation step in the original DRS algorithm. Additionally, we apply a statistical test to a set of $4$-dimensional vectors generated using the original DRS algorithm. We demonstrate that these vectors are not necessarily drawn from a uniform distribution.

We extend the DRS algorithm by introducing a method called the Dirichlet-Rescale-Constraints (DRSC) algorithm to allow for the generation of vectors with fixed sum that satisfy both linear and nonlinear constraints. We also resolve the issue of non-uniformity, such that vectors generated by the DRSC algorithm follow a uniform distribution.

The remainder of this paper is organised as follows. In \Cref{section2:concepts}, we explain mathematical concepts needed for the DRS and DRSC algorithms. We present a simplified DRS algorithm in \Cref{section3:DRS} and analyse its correctness in \Cref{section4:analysingCorrectness}. We introduce the DRSC algorithm and perform computational experiments in Sections~\ref{section5:DRSC} and \ref{section6:computationalResults}, respectively. Finally, in \Cref{section7:conclusion} we conclude and discuss the implications of our findings for the use of the DRS algorithm.

\section{Mathematical concepts}\label{section2:concepts}
\begin{figure}[h]
    \centering
    \includegraphics[width=0.55\textwidth]{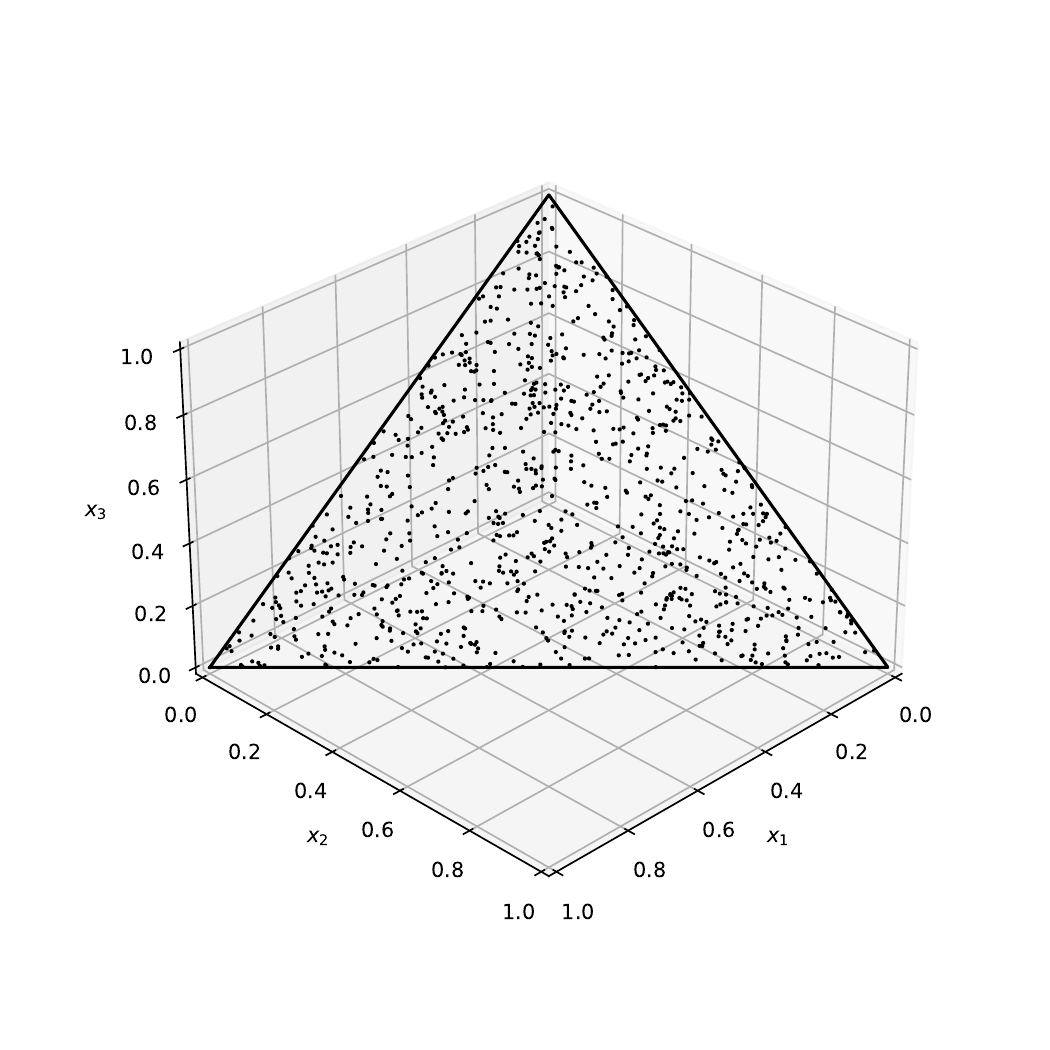}
    \caption{A sample drawn from a standard simplex. Each point corresponds to a vector.}
    \label{fig:standardsimplex}
\end{figure}

A key concept in this paper is the notion of a simplex $S$, which is a generalisation of a triangle in higher dimensions. {Specifically, an $n-1$-simplex is defined as the convex hull of $n$ affinely independent $n$-dimensional vertices.} We focus on two types of simplices: regular and standard. A simplex is called regular if all the edges have the same length. The standard simplex is a regular simplex with standard unit vectors, as illustrated in \Cref{fig:standardsimplex}. All simplices are affine, meaning any simplex can be transformed into another simplex of the same dimensionality through an affine transformation. This property is useful, because a uniformly distributed sample retains its uniform distribution after applying an affine transformation. 

The Dirichlet distribution, parametrised by a parameter vector $\boldsymbol{\alpha}$, generalises the Beta distribution \citep{johnson1972}. The flat Dirichlet distribution, obtained when $\boldsymbol{\alpha}=\matr{1}$, is equivalent to a uniform distribution over the standard simplex. {Thus, sampling from a flat Dirichlet distribution yields a vector that sums to $1$ and that is drawn uniformly from the simplex \citep{griffin2020}. Each point in \Cref{fig:standardsimplex} represents a vector that is randomly generated from such a distribution.}

\FloatBarrier

\section{The Dirichlet-Rescale algorithm}\label{section3:DRS}
In \Cref{alg:DRS}, we present a simplified version of the Dirichlet-Rescale (DRS) algorithm, as proposed by \cite{griffin2020}, to generate vectors $\matr{x}$ with non-negative entries that sum to~1 and satisfy lower and upper bounds on the entries. Vector $\matr{x}$ that do not sum to 1 can be transformed accordingly \citep{griffin2020}. The input to the DRS algorithm is a set of constraints $J$, which specifies individual lower and upper bounds for each entry.

When a constraint is violated, a restart can always be performed by drawing a new sample from the flat Dirichlet distribution, similar to the principles of the acceptance-rejection method. However, this approach is unlikely to be efficient in practice. To improve efficiency, the DRS algorithm applies affine transformations to project the vectors towards the feasible region until all constraints are satisfied.
\begin{algorithm}[h]
\caption{The Dirichlet-Rescale algorithm}
\label{alg:DRS}
\begin{algorithmic}[1]
\State Generate a vector $\matr{x}$ on the standard simplex $S$ by sampling it from the flat Dirichlet distribution. 
\If{all constraints in $J$ are satisfied by $\matr{x}$}
\State Stop. 
\EndIf
\State Construct an induced simplex $S'$, based on the violated constraints in $J$.  
\State An affine transformation maps $\matr{x}$ from the induced simplex $S'$ onto the standard simplex $S$. Go to line 2.
\end{algorithmic}
\end{algorithm}
\FloatBarrier

\begin{figure}[h]
    \centering
    \begin{subfigure}[b]{0.44\textwidth}
        \centering
        \includegraphics[width=\textwidth]{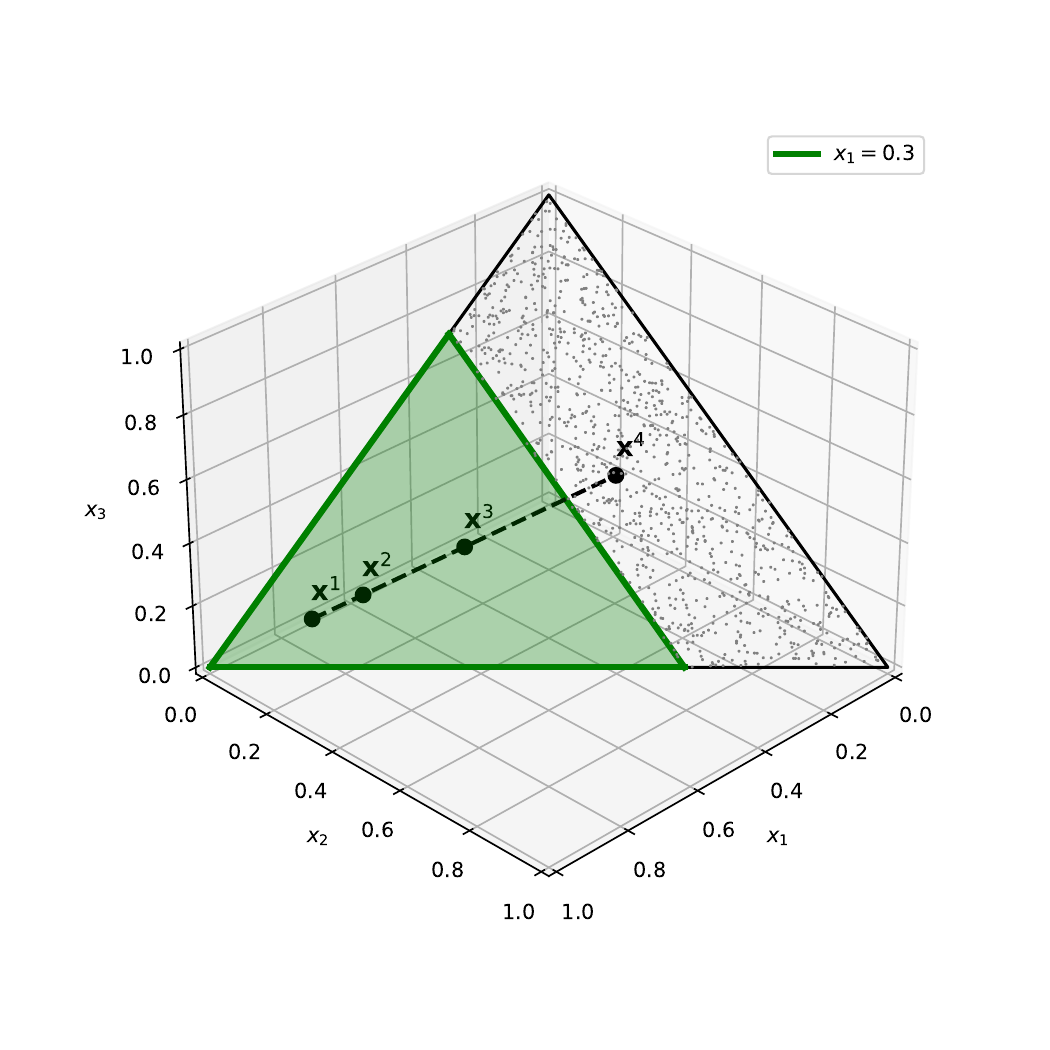}
        \caption[]%
        {{\small A standard simplex with constraint $x_1\leq 0.3$. Starting from $\matr{x}^1$, repeatedly applying an affine transformation mapping the green simplex onto the standard simplex yields $\matr{x}^2$, $\matr{x}^3$ and $\matr{x}^4$.}}    
        \label{fig:upperboundDRS}
    \end{subfigure}
    \hfill
    \begin{subfigure}[b]{0.44\textwidth}  
        \centering 
        \includegraphics[width=\textwidth]{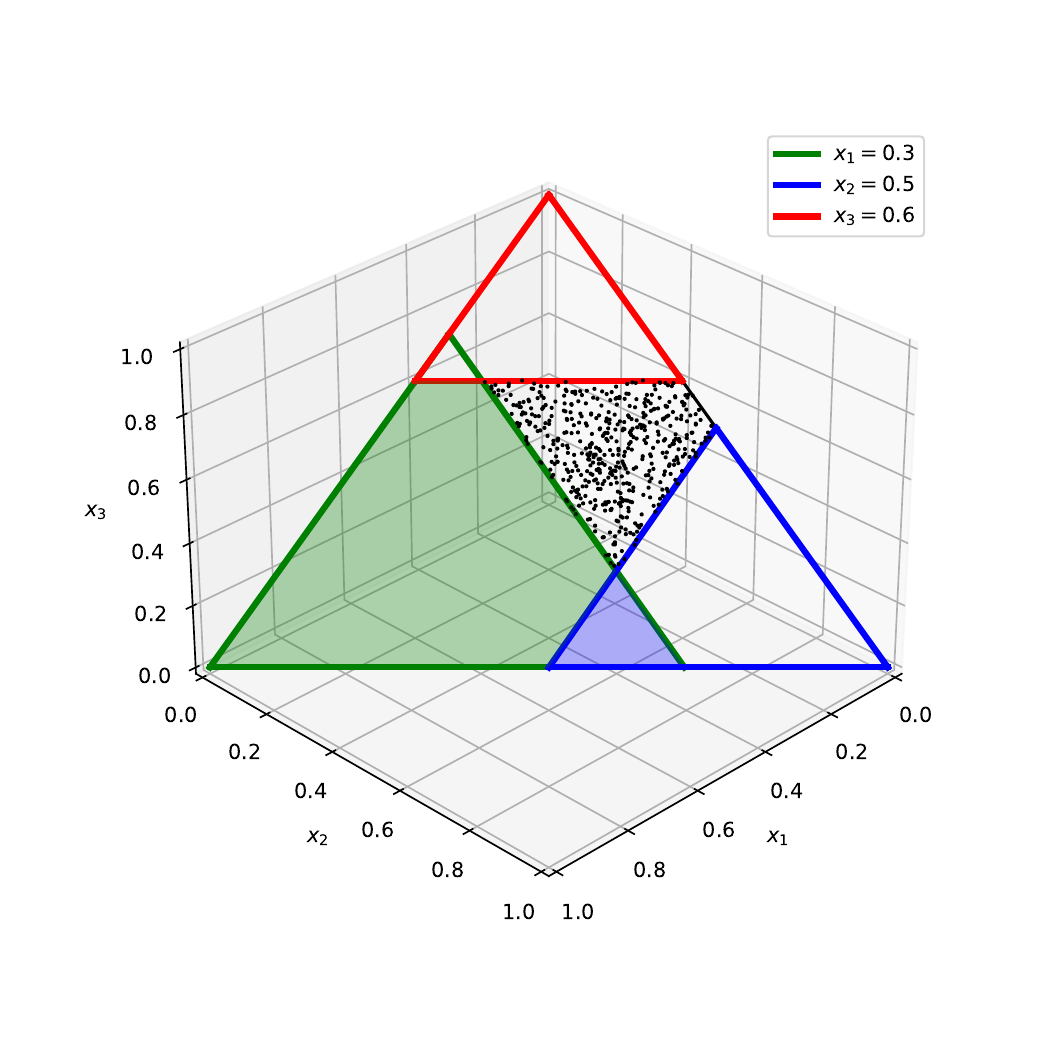}
        \caption[]%
        {{\small A standard simplex with constraints $x_1\leq 0.3$, $x_2\leq 0.5$ and $x_3\leq 0.6$. The overlapping constraints $x_1\leq 0.3$ and $x_2\leq 0.5$ induce the small blue shaded simplex. However, constraint $x_2\leq 0.5$ no longer induces a simplex, which is green shaded.}}    
        \label{fig:standardsimplex_withcrossingconstraints}
    \end{subfigure}
    \caption{Standard simplices with lower and upper bound constraints.}
    \label{fig:DRS}
\end{figure}

{We illustrate the DRS algorithm using \Cref{fig:upperboundDRS}, which considers a $3$-dimensional standard simplex $S$ with constraint $x_1\leq 0.3$. Initially, vectors are drawn from a uniform distribution over the standard simplex $S$, which is equivalent to sampling from a flat Dirichlet distribution. Suppose that we draw $\matr{x}^1$, which violates the constraint $x_1\leq 0.3$. Based on the violated constraint, we create an induced simplex $S'=\{x\in S:   x_1\geq 0.3 \}$. In line 5 we construct the simplex $S'$ induced by the constraint $x_1\leq 0.3$, which is marked in green in \Cref{fig:upperboundDRS}. By construction, the vector $\matr{x}^1$ lies within simplex $S'$. Since $S'$ is a subset of $S$, vectors within $S'$ are also a sample from a uniform distribution. In line 6, an affine transformation maps $\matr{x}^1$ from simplex $S'$ to the standard simplex $S$, preserving the uniform distribution. This transformation projects vector $\matr{x}^1$ towards the feasible region. As depicted in \Cref{fig:upperboundDRS}, this process maps $\matr{x}^1$ to $\matr{x}^2$. These steps are repeated until all constraints are satisfied, which happens when $\matr{x}^4$ is reached.}

\FloatBarrier

The aim of the DRS algorithm is to draw vectors uniformly from the set of possible vectors constrained by the individual bounds in $J$, meaning that each vector within the feasible region is equally likely to be selected. However, it turns out that this uniformity property stated by \cite{griffin2020} does not hold. In the next section, we analyse scenarios where the DRS algorithm does not necessarily exhibit the uniform property.

\section{On the uniformity of draws from the Dirichlet-Rescale algorithm}\label{section4:analysingCorrectness}
{When multiple constraints are violated, it may be impossible to partition the infeasible region into simplices. In that case, the DRS algorithm proposed by \cite{griffin2020} does not necessarily provide vectors drawn from a uniform distribution.} For instance, consider a vector that violates two constraints: $x_1\leq 0.3$ and $x_2\leq 0.5$, as shown in \Cref{fig:standardsimplex_withcrossingconstraints}. Here, the blue shaded area induces a simplex, but this causes the green shaded area to no longer induce a simplex. Consequently, the DRS algorithm is no longer guaranteed to generate a uniform distribution of vectors.

We provide two procedures, namely an analytical procedure and a statistical test, to demonstrate that the DRS algorithm as proposed by \cite{griffin2020} may not always generate vectors that are drawn from a uniform distribution. First, we provide an analytical procedure showing that $3$-dimensional vectors generated using the simplified DRS algorithm described in \Cref{alg:DRS} do not, in general, follow from a uniform distribution.

We note that the original DRS algorithm proposed by \cite{griffin2020} contains a transformation step, which switches the roles of the standard simplex with the simplex induced by all the lower and upper bound constraints. This transformation step may have the unintentional effect that $3$-dimensional vectors appear to be drawn from a uniform distribution. In a second procedure, we apply a statistical test to show that the original DRS algorithm from \cite{griffin2020} does not always generate $4$-dimensional vectors from a uniform distribution.

\subsection{Analytical procedure: tiling the simplex}\label{appendixA:tiling}
The simplified DRS algorithm as outlined in \Cref{alg:DRS} is deterministic. Given a starting vector $\matr{x}^1$ we can repeatedly apply the affine transformation steps from the DRS algorithm to determine in how many steps the vector arrives in the feasible region, as illustrated in \Cref{fig:upperboundDRS}.

Consider a standard simplex $S$ with mass $m(S)$ and an upper bound vector $\matr{u}$, where there is a possibility of violating two constraints simultaneously. We perform the DRS algorithm in ``reverse'' by partitioning the standard simplex $S$ into multiple regions, that we refer to as tiles. Specifically, we partition the simplex $S$ using tiles $T_i$, such that all vectors within a tile $T_i$ reach the feasible region in exactly $i$ steps. The feasible region, induced by the constraint $\matr{u}$, is denoted by $T_0$. With this procedure, we demonstrate that applying the DRS algorithm may lead to a non-uniform distribution of projected mass $m(S)$ over the feasible region $T_0$. This implies that the vectors in $T_0$ are not drawn from a uniform distribution.

\begin{figure}[h]
    \centering
    \includegraphics[width=0.67\textwidth]{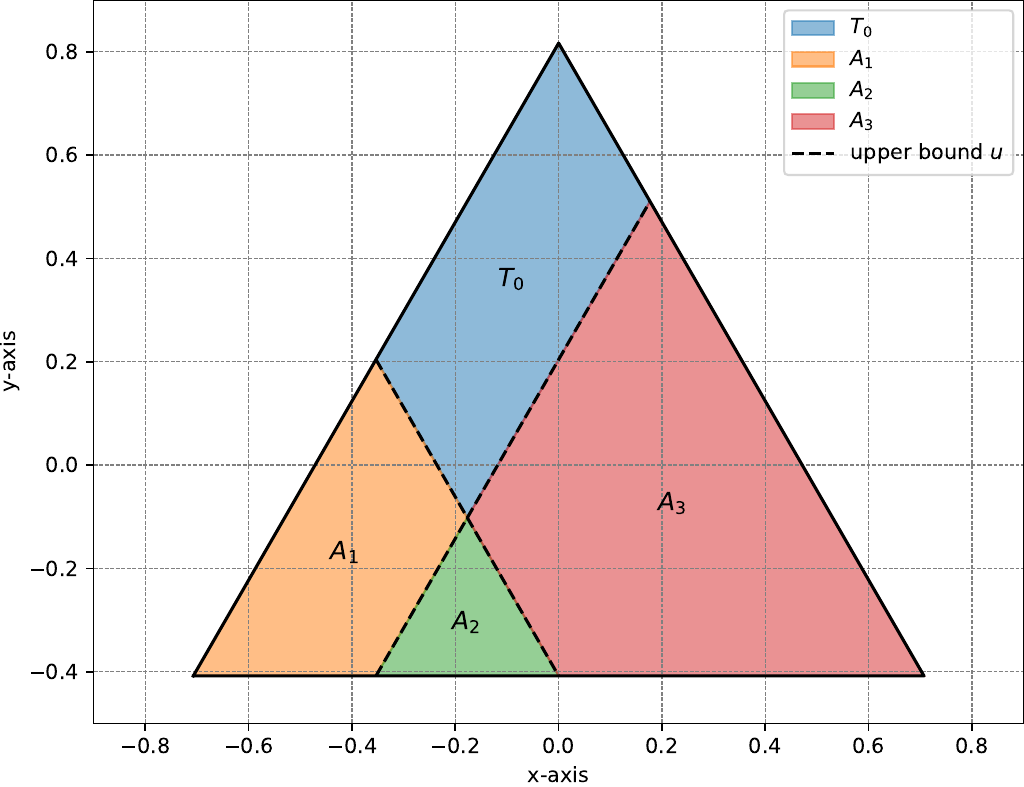}
    \caption{The projection of a standard simplex with an upper bound vector $\matr{u}$, resulting in two constraints. The constraints induce the feasible region $T_0$. Vectors generated outside the feasible region are transformed using $A_1$, $A_2$ or $A_3$.}
    \label{fig:actionSpace}
\end{figure}

To illustrate the procedure, consider a $3$-dimensional setting with an upper bound vector $u=[0.5, 0.25, 1]$ on $\matr{x}$. For ease of interpretation and without loss of generality, we project the $3$-dimensional simplex onto a $2$-dimensional space using an orthogonal transformation. \Cref{fig:actionSpace} illustrates the feasible region $T_0$ and three areas induced by the constraints. Let $\mathcal{A}$ be the set of possible transformations in a DRS step. In this example, each infeasible region corresponds with one possible transformation, namely $A_1$, $A_2$ or $A_3$. For instance, when a vector $\matr{x}$ lies in the infeasible region associated with $A_1$, we apply transformation $A_1$ to $\matr{x}$.

\FloatBarrier

Starting from $T_0$, we recursively apply $i$ ``reverse'' DRS steps, while taking into account the bound $\matr{u}$, to obtain tiles $T_i^{(A_1, \dots, A_i)}$, where $A_j$ denotes the $j$-th transformation $A\in \mathcal{A}$. Note that the same transformation may be applied more than once.
The notation $T_i^{(A_1, \dots, A_i)}$ indicates that applying the sequence of transformations $A_1, \dots, A_i$ to tile $T_i^{(A_1, \dots, A_i)}$ projects all vectors onto the projected tile $P(T_i^{(A_1, \dots, A_i)})$ in exactly $i$ steps. The projected tile $P(T_i^{(A_1, \dots, A_i)})$ is by definition either equal to or contained within $T_0$. For ease of notation, we leave out the transformations and simply refer to tiles $T_i$ and their projections $P(T_i)$, even though these tiles may have a different shape and mass.

For each tile $T_i$, we compute the mass $m(T_i)$, which in a 3-dimensional simplex corresponds to the area of the tile. Since we expect some projected tiles $P(T_i)$ to be a strict subset of $T_0$, we partition $T_0$ in regions such that within a region the mass is equally distributed. Using the mass of the standard simplex, $S$, we analytically compute a target mass per region, based on the assumption of uniformity. In addition, based on the mass of the projected tiles, $P(T_i)$, we determine analytically the realised mass per region. If the DRS algorithm generates a uniform distribution of vectors, the target and realised mass should be equal in every region. The tiling of the simplex, and the analytical computation of the target and realised mass per region are summarised in \Cref{alg:fractal}. 

\begin{algorithm}[h]
\caption{Tiling the simplex}
\label{alg:fractal}
\begin{flushleft}
    \textbf{Inputs:} A standard simplex $S$ with mass $m(S)$, a set of transformation matrices $\mathcal{A}$ and an upper bound vector $\matr{u}$\\
\end{flushleft}
\begin{algorithmic}[1]
\State Find tiles $T_i^{(A_1, \dots, A_i)}$ by recursively applying $i$ ``reverse'' DRS steps using transformations in $\mathcal{A}$, while taking into account the bound $\matr{u}$. 
\State Project the tiles $T_i^{(A_1, \dots, A_i)}$ onto $P(T_i^{(A_1, \dots, A_i)})$.
\State Partition $T_0$ in smaller regions and determine the target mass per region. 
\State Based on the projected tiles $P(T_i^{(A_1, \dots, A_i)})$ calculate the realised mass per region
\end{algorithmic}
\end{algorithm}
\FloatBarrier

\begin{figure}[h]
    \centering
    \includegraphics[width=0.67\textwidth]{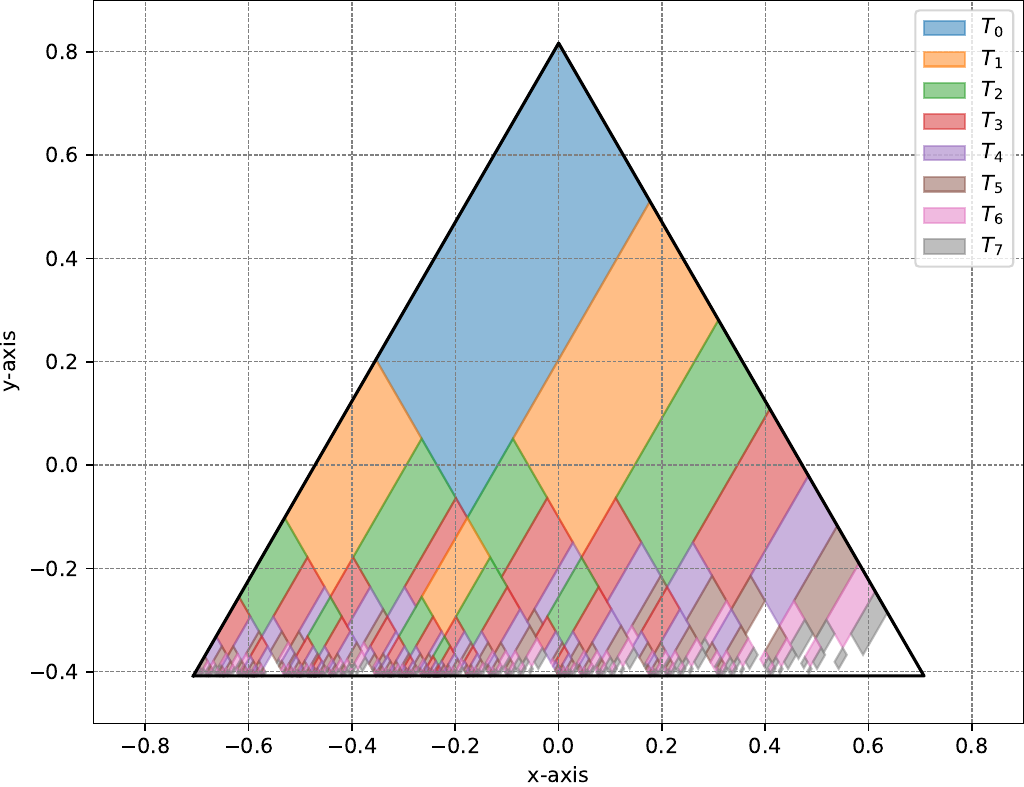}
    \caption{Tiling of the simplex. Each tile $T_i$ indicates that $i$ DRS steps must be performed to reach the feasible region $T_0$.}
    \label{fig:fractal_tiling}
\end{figure}

By applying \Cref{alg:fractal} with a maximum of $i=7$ DRS steps we obtain the tiling shown in \Cref{fig:fractal_tiling}. Note that the shapes of some projected tiles $T_i$ differ from tile $T_0$, corresponding to the feasible region. 

\FloatBarrier

\begin{figure}[h]
    \centering
    \includegraphics[width=0.67\textwidth]{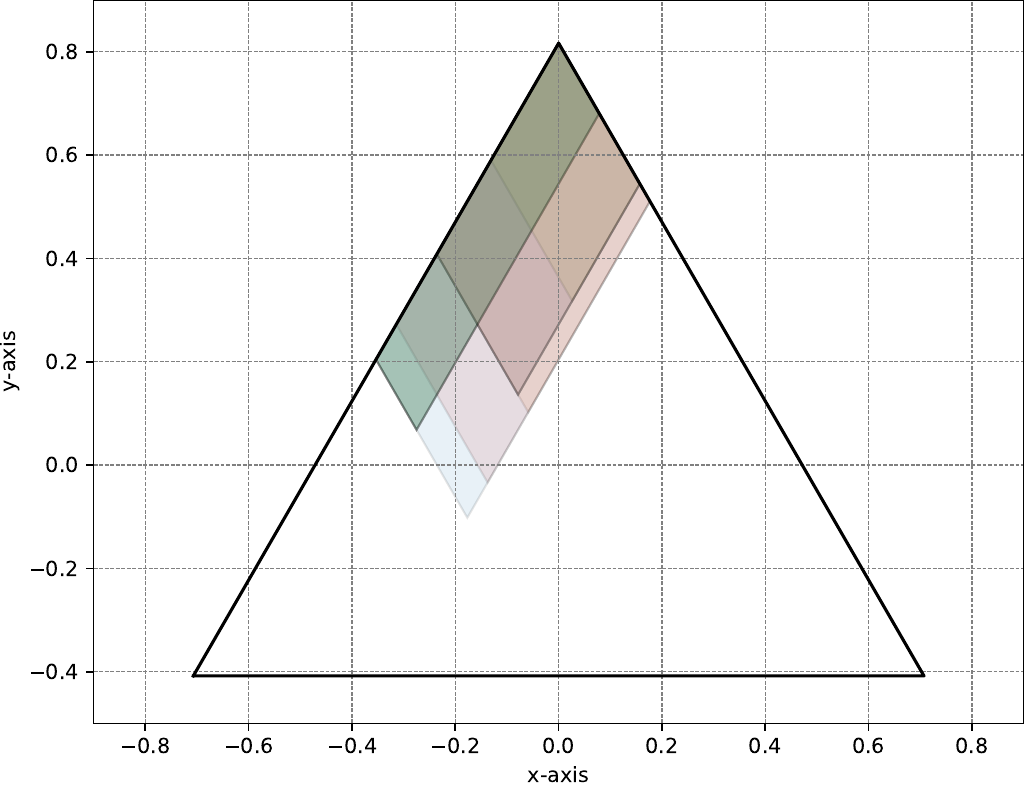}
    \caption{Projected tiles $P(T_i$) with a unique shape and the feasible region $T_0$.}
    \label{fig:fractal_overlap}
\end{figure}

\Cref{fig:fractal_overlap} illustrates the projected tiles $P(T_i)$ with a unique shape. The shaded regions seem to indicate that fewer vectors are projected into the bottom corner of tile $T_0$ compared to its top corner.

\FloatBarrier

\begin{figure}[h]
    \centering
    \includegraphics[width=0.67\textwidth]{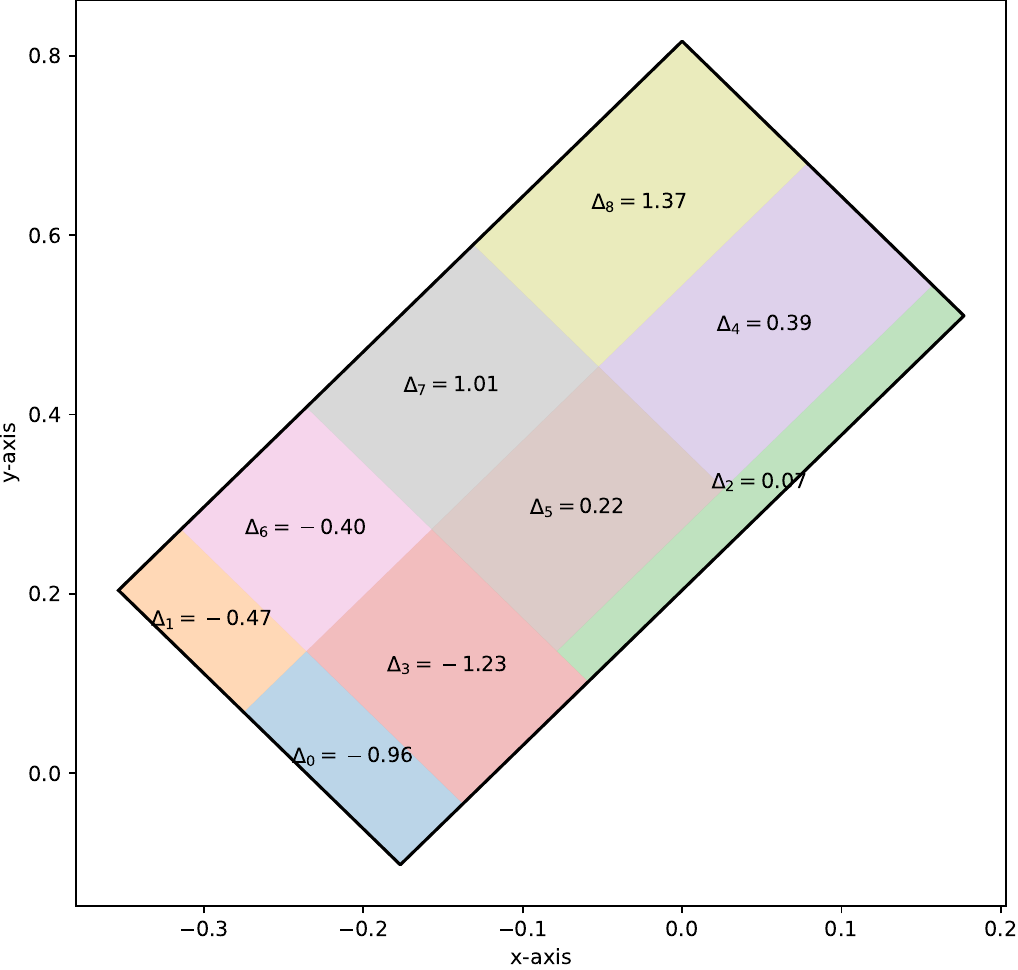}
    \caption{The analytically determined differences ($\Delta$) between the realised and target mass in percentages in the feasible region $T_0$.}
    \label{fig:fractal_delta}
\end{figure}

Based on \Cref{fig:fractal_overlap}, we can partition $T_0$ into nine regions, such that within a region the analytically determined mass remains uniformly distributed. This partitioning is illustrated in \Cref{fig:fractal_delta}. For each region, we compute a value $\Delta$, representing the percentage difference between the realised and target mass after having performed $i=7$ DRS steps. In the bottom corner the realised mass is approximately 1 percentage point less than the target, whereas in the top corner the realised mass exceeds the target by over 1 percentage point.

\Cref{tab:fractal} contains detailed results on the {realised and target } mass in percentages. The realised and target percentages both account for 96.52\% of the total mass. The residual of 3.48\% is unallocated because the corresponding mass requires more than $i=7$ DRS steps to reach the feasible region $T_0$. Performing additional DRS steps would allocate the residual mass, leading to more accurate realised and target mass percentages. However, even in an optimal scenario where the residual of 3.48\% is fully allocated, this would not be enough to compensate for the sum of absolute differences, which totals 6.12\%.

To conclude, when two or more violated constraints can be violated, the simplified DRS algorithm described in \Cref{alg:DRS} does not necessarily generate vectors from a uniform distribution. 

\begin{table}[h]
\caption{Realised and target mass in percentages per region. The residual mass corresponds to all tiles $T_i$ with $i> 7$.}
    \label{tab:fractal}
    \centering
\begin{tabular}{lrrr}
\toprule
region                    & realised (in \%) & target (in \%)& difference $\Delta$ (in \%)\\ \midrule
0                         & 5.00     & 5.96     & -0.96      \\
1                         & 4.29     & 4.77     & -0.47      \\
2                         & 7.22     & 7.15     & 0.07       \\
3                         & 10.68    & 11.92    & -1.23      \\
4                         & 16.28    & 15.89    & 0.39       \\
5                         & 12.93    & 12.71    & 0.22       \\
6                         & 9.14     & 9.53     & -0.40      \\
7                         & 13.72    & 12.71    & 1.01       \\
8                         & 17.26    & 15.89    & 1.37       \\ \midrule
total (in regions 0 to 8) & 96.52 & 96.52 & \\ \midrule
residual & 3.48 & 3.48 &  \\ \midrule
sum of absolute differences &          &          & 6.12      \\ \bottomrule
\end{tabular}
\end{table}
\FloatBarrier

\subsection{Goodness of fit test}\label{appendixB:goodnessoffittest}
In this section, we analyse the original DRS algorithm as proposed by \cite{griffin2020}. Specifically, we perform a $\chi^2$ goodness of fit test to evaluate the hypothesis that vectors from the DRS algorithm are generated from a uniform distribution in an $n$-dimensional space. 

\cite{griffin2020} conducted a general statistical test to assess uniformity, while we perform a more specific test. Their approach involves a two-sample test, where they use a method proven to follow a uniform distribution to construct a reference sample. Their analysis is applied on $10$-dimensional vectors with a randomly drawn upper bound constraint. They found no evidence to reject the null hypothesis that the samples were similar, suggesting that the DRS algorithm generates a uniform distribution of vectors. In contrast, because we have a clear understanding of the cause of the non-uniformity, we are able to perform a more specific statistical test.

The main idea of the statistical test that we perform, is to partition the feasible region into equally sized bins and evaluate whether the vectors are distributed evenly across these bins. In our statistical test, we consider $4$-dimensional vectors $\matr{x}$ drawn using the original DRS algorithm. Let $u=[0.1, 0.5, 0.8, 1]$ be the upper bound vector on $\matr{x}$. We draw two samples, of size $s_1=10,000$ and $s_2=10,000,000$, respectively. 

Samples $s_1$ and $s_2$ are bijectively projected onto an $n-1$-dimensional space by dropping the $n$-th coordinate. Using the first sample $s_1$, we find the convex hull $C$ (by using a Delaunay triangulation). The space $[0, 1]^{n-1}$ is divided into $n_b=25$ equally sized bins for each dimension, leading to $n_b^{n-1}$ potential bins. 

We iterate through all possible bins, retaining only those that are entirely contained within the convex hull $C$, resulting in the set $B_{}$. After, we count the number of vectors from the second sample, $s_2$, that fall into each bin $b\in B_{}$, yielding observed counts $O_b$. The total number of vectors considered in the $\chi^2$ test is given by $\sum_{b\in B_{}} O_b$. Since we have equally sized bins, the expected number of vectors per bin is calculated as $E = \frac{\sum_{b\in B_{}} O_b}{|B_{}|}$. 

The $\chi^2$ statistic can be computed on sample $s_2$ as
\begin{align*}
    \chi^2 = \sum_{b\in B_{}} \frac{(O_b - E)^2}{E}
\end{align*}
with $|B|-1$ degrees of freedom.

In our experiment, we obtain $|B|=178$, $\chi^2=214.8879$ and a $p$-value of 0.0274. Using a significance level of $5\%$, we reject the null hypothesis. We conclude that the vectors from sample $s_2$ obtained using the DRS algorithm are not from a uniform distribution.

\section{The Dirichlet-Rescale-Constraints algorithm}\label{section5:DRSC}
{In this section, we propose an algorithm that we call Dirichlet-Rescale-Constraints (DRSC). The DRSC algorithm is an extension of the DRS algorithm designed for sampling vectors with a fixed sum that handles linear and nonlinear constraints.} The DRSC algorithm follows a similar structure as \Cref{alg:DRS}, however, set $J$ may now include linear and nonlinear constraints.

\begin{figure}[h]
    \centering
     \begin{subfigure}[b]{0.44\textwidth}  
        \centering 
        \includegraphics[width=\textwidth]{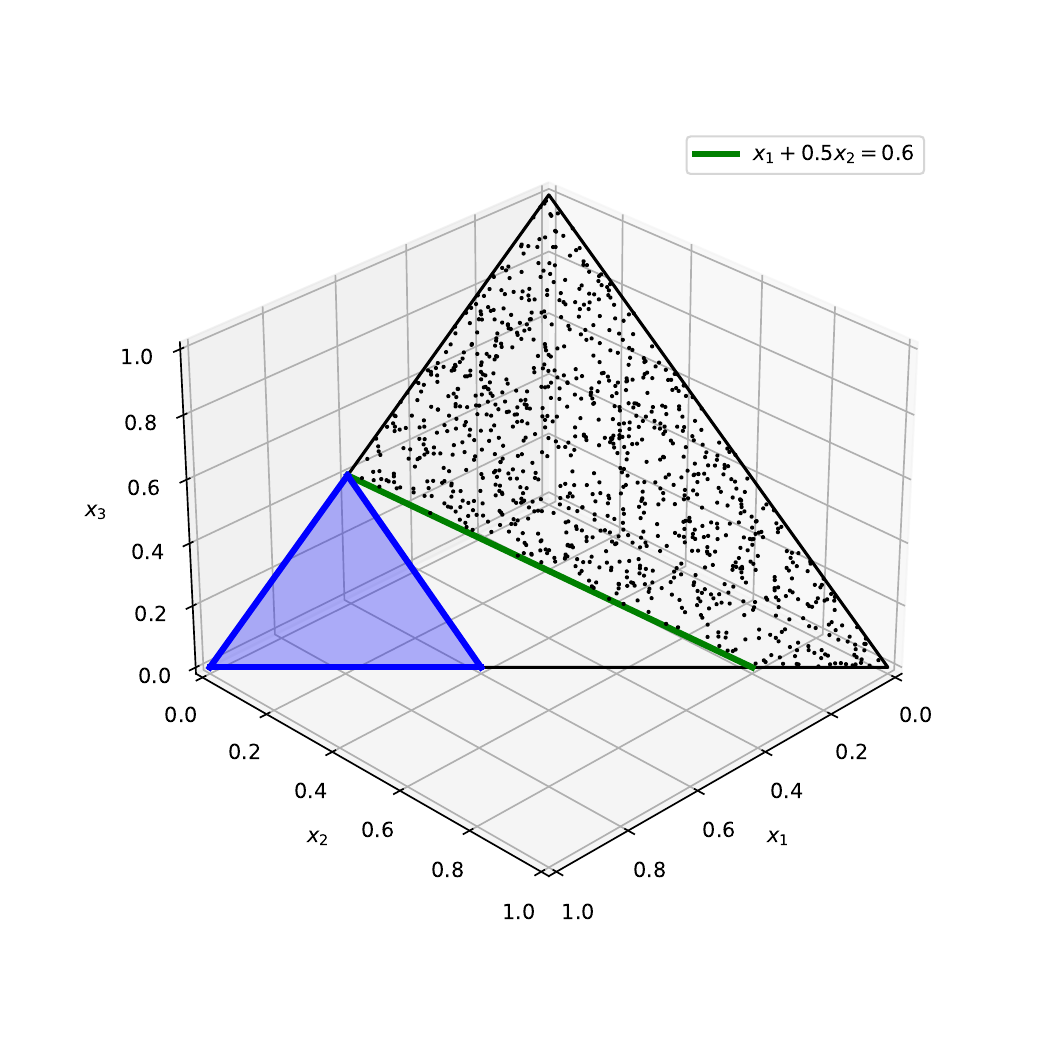}
        \caption[]%
        {{\small A standard simplex with inequality $x_1+0.5x_2\leq 0.6$. The blue simplex is the largest possible regular simplex induced by the constraint.}}    
        \label{fig:standardsimplex_withinequalityconstraints}
    \end{subfigure}   
    \hfill
    \begin{subfigure}[b]{0.44\textwidth}
        \centering
        \includegraphics[width=\textwidth]{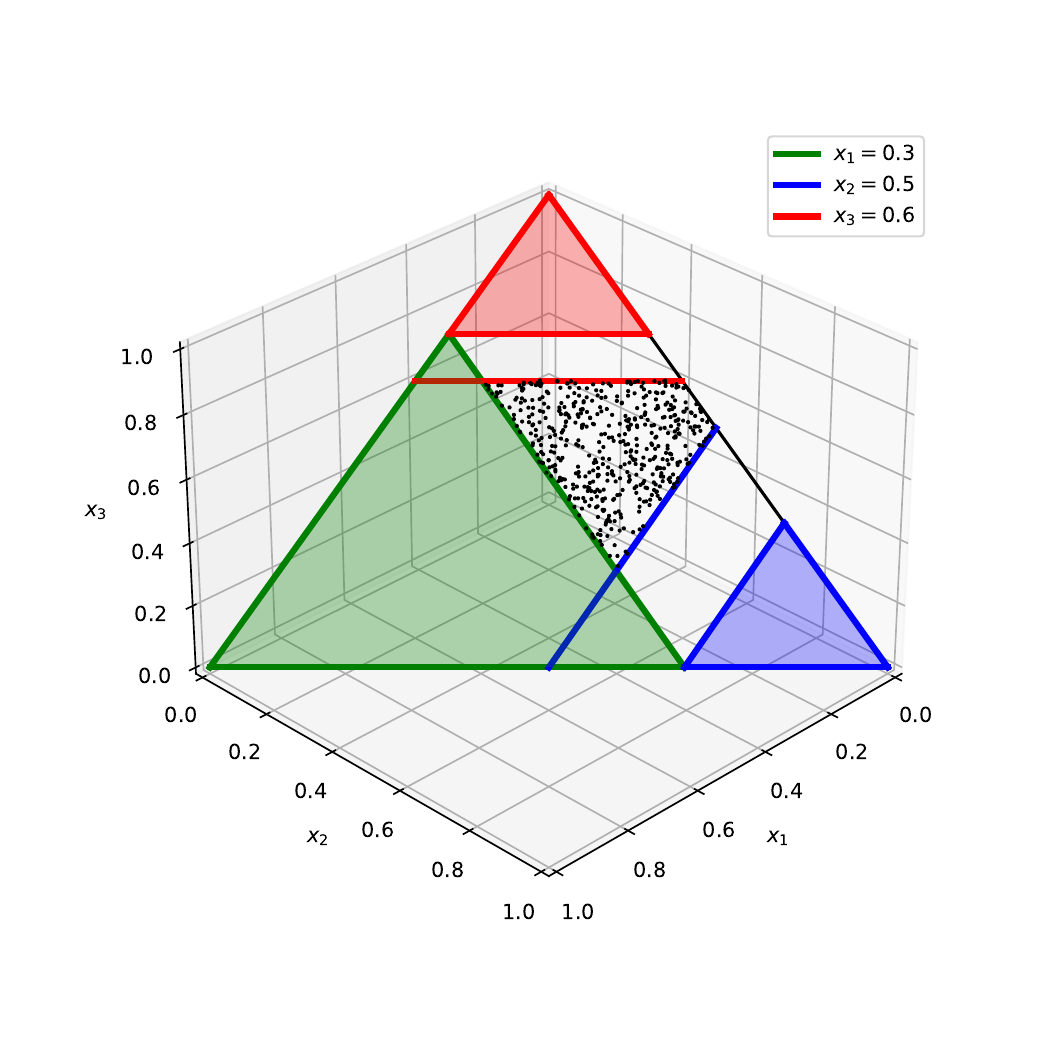}
        \caption[]%
        {{\small A standard simplex with constraints $x_1\leq 0.3$, $x_2\leq 0.5$ and $x_3\leq 0.6$. The green, blue and red shaded areas correspond to non-overlapping regular simplices induced by the constraints.}}    
    \label{fig:standardsimplex_withcrossingconstraints_largestSimplex}
    \end{subfigure}
    \caption{Standard simplex with induced regular simplices.}
    \label{fig:DRSC_linear}
\end{figure}

Let $J_{lin}\subseteq J$ be the set of linear constraints. Similar to the DRS algorithm, when a vector violates a constraint from $J_{lin}$, we construct an induced simplex and apply an affine transformation to project the vector onto the standard simplex. The DRSC algorithm addresses two key challenges. The first challenge, as illustrated in \Cref{fig:standardsimplex_withcrossingconstraints}, arises when there are overlapping constraints, since this may lead to a non-uniform distribution of generated vectors. The second challenge involves handling any linear inequality, as shown in \Cref{fig:standardsimplex_withinequalityconstraints}, where the inequality $x_1+0.5x_2\leq 0.6$ induces an irregular simplex.

Our main idea to address the first challenge is illustrated in \Cref{fig:standardsimplex_withcrossingconstraints_largestSimplex}. We resolve the issue arising from overlapping constraints by identifying non-overlapping regular simplices. To address the second challenge of handling general linear inequalities, we propose finding the largest possible regular simplex induced by the linear constraints. For instance, in \Cref{fig:standardsimplex_withinequalityconstraints}, the largest induced regular simplex is shaded in blue. To solve both key challenges, we aim to find large non-overlapping regular simplices.

When a generated vector lies within an induced regular simplex, we can apply an affine transformation mapping the induced simplex onto the standard simplex. However, when the generated vector violates a constraint and does not lie within any induced simplex, we restart the algorithm by drawing from the flat Dirichlet distribution again.

\FloatBarrier

Linear inequality constraints can be expressed in vector notation as $\matr{a}^\top \matr{x}\leq {b}$. Equality constraints, $\matr{a}^\top \matr{x} = b$, can similarly be modelled as inequalities by reformulating them as
\begin{align*}
   \matr{a}^\top \matr{x} &\leq b,\\
    - \matr{a}^\top \matr{x} &\leq b.
\end{align*}

Consider a linear constraint of the form $\matr{a}_j^\top \matr{x} \leq b_j$ for $j\in J_{lin}$. Our goal is to find large non-overlapping simplices that are induced by all of the linear constraints in $J_{lin}$. Specifically, for each dimension $i=1,\dots, n$ we seek to maximise a constant $\theta_i$ that induces a simplex $S_i$ through the constraint $x_i\geq \theta_i$. When a vector $\matr{x}$ lies within the induced simplex $S_i$, it violates at least one constraint from $J_{lin}$.

For each dimension $i=1,\dots, n$, we solve a linear programming problem to obtain a value $\theta_i$. These values $\theta_i$ serve as inputs for the subsequent linear programs, ensuring that the simplices are non-overlapping.

Let $\matr{e}_i$ denote the unit vector whose $i$-th entry is 1, and all other entries are 0. Furthermore, let $\matr{x}^k$ for $k=1,\dots, i$ be a vector of variables, specifically, $x_k^k$ denotes the $k$-th element of vector $\matr{x}^k$. The largest regular induced simplex $S_i$, corresponding to $x_i\geq \theta_i$, can be identified by solving the following linear programming formulation
\begin{align}
    \theta_{i}=\max\ & \max_{k=1,\dots, i} \matr{e}_i^\top\matr{x}^k, \label{form:DRSC_1} \\
    \text{s.t. }    & \matr{1}^\top\matr{x}^i = 1,  \label{form:DRSC_2} \\
    & \matr{0} \leq \matr{x}^i \leq \matr{1},  \label{form:DRSC_3}\\
 &   \matr{a}_j^\top \matr{x}^i \leq  b_j, && \forall j\in J_{lin}, \label{form:DRSC_4} \\
   & \matr{1}^\top\matr{x}^k = 1, && \forall k=1,\dots, i-1, \label{form:DRSC_5} \\
    & \matr{0} \leq \matr{x}^k \leq \matr{1},&& \forall k=1,\dots, i-1.  \label{form:DRSC_6}\\
        & \theta_k \leq {x}_k^k ,&& \forall k=1,\dots, i-1.  \label{form:DRSC_7}
\end{align}
We maximise (\ref{form:DRSC_1}) to determine the largest possible value for $\theta_i$. The feasible region of $\matr{x}^i$ is modelled by (\ref{form:DRSC_2})-(\ref{form:DRSC_4}). Specifically, (\ref{form:DRSC_2}) and (\ref{form:DRSC_3}) model a standard simplex, while (\ref{form:DRSC_4}) enforces that the linear constraints $j\in J_{lin}$ are satisfied. The existing induced simplices $S_k$ with $k=1,\dots, i-1$ are modelled by (\ref{form:DRSC_5})-(\ref{form:DRSC_7}) using variables $\matr{x}^k$. The standard simplex is modelled by (\ref{form:DRSC_5}) and (\ref{form:DRSC_6}), while (\ref{form:DRSC_7}) ensures that the induced simplex $S_i$ does not overlap with the previously constructed simplex $S_k$.

Note that the objective (\ref{form:DRSC_1}) can be linearised and the formulation (\ref{form:DRSC_1})-(\ref{form:DRSC_7}) can be decomposed, as the $\matr{x}_k$ variables are independent from each other. In practice, for each vector of variables $\matr{x}^k$ a separate linear programming formulation is solved. Furthermore, note that the described approach is a heuristic approach, since it optimizes the induced simplex in a fixed order.

In the implementation it is not possible to exactly satisfy an equality constraint, as the probability of satisfying an equality constraint would be 0. In addition, it may happen that multiple inequalities induce an equality constraint. Therefore, we model inequality constraints and equality constraints as 
\begin{align*}
    \matr{a}^\top \matr{x}&\leq {b} + \epsilon_{ineq},\\
       \matr{a}^\top \matr{x} &\leq b + \epsilon_{eq},\\
    - \matr{a}^\top \matr{x} &\leq b+ \epsilon_{eq},
\end{align*}
where $\epsilon_{ineq}$ and $\epsilon_{eq}$ denote the feasibility tolerance for inequality and equality constraints, respectively.

Let us now consider the set of nonlinear constraints $J_{non}\subseteq J$. When a nonlinear constraint $f_j(\matr{x})\leq b$ is violated, the DRSC algorithm restarts by drawing from a flat Dirichlet distribution. The DRSC algorithm is summarised in \Cref{alg:DRSC}, with a similar structure as \Cref{alg:DRS}.

\begin{algorithm}[h]
\caption{The Dirichlet-Rescale-Constraints algorithm}
\label{alg:DRSC}
\begin{algorithmic}[1]
\State Find the induced simplices $S_i$ for each dimension $i=1,\dots, n$ using (\ref{form:DRSC_1})-(\ref{form:DRSC_5}).
\State Generate a vector $\matr{x}$ on the standard simplex $S$ by sampling it from the flat Dirichlet distribution. 
\If{all constraints in $J$ are satisfied by $\matr{x}$}
\State Stop. 
\ElsIf{the vector $\matr{x}$ is inside an induced simplex $S_i$}
\State An affine transformation maps $\matr{x}$ from the induced simplex $S_i$ onto the standard simplex~$S$. Go to line 3.
\EndIf
\State Go to line 2.

\end{algorithmic}
\end{algorithm}
\FloatBarrier

The vectors generated by the DRSC algorithm are drawn from a uniform distribution, as formalised in Property~\ref{property:uniformDRSC}.
\begin{property}\label{property:uniformDRSC}
    The vectors generated by the DRSC algorithm are draws from a uniform distribution over the region defined by the constraint set $J$.
\end{property}

To preserve Property~\ref{property:uniformDRSC}, it is important to first verify whether $\matr{x}$ lies within an induced simplex before rejecting it based on a nonlinear constraint. Reversing this order would disrupt the induced simplex, such that the uniformity property no longer holds. Thus, line 8 is reached either when vector $\matr{x}$ violates a constraint in $J_{lin}$ but does not lie within the region of any induced simplex $S_i$, or when $\matr{x}$ violates a nonlinear constraint.

\section{Computational results}\label{section6:computationalResults}
In this section, we present the numerical results on generating contingency tables. We demonstrate that the Dirichlet-Rescale-Constraints algorithm is able to generate vectors with fixed sum while having linear and nonlinear constraints. The computational experiments are performed on a 1.8 GHz Intel Core i7 processor with 16.0 GB RAM.

\begin{figure*}[h]
    \centering
    
    \begin{subfigure}[b]{0.44\textwidth}  
        \centering 
        \includegraphics[width=\textwidth]{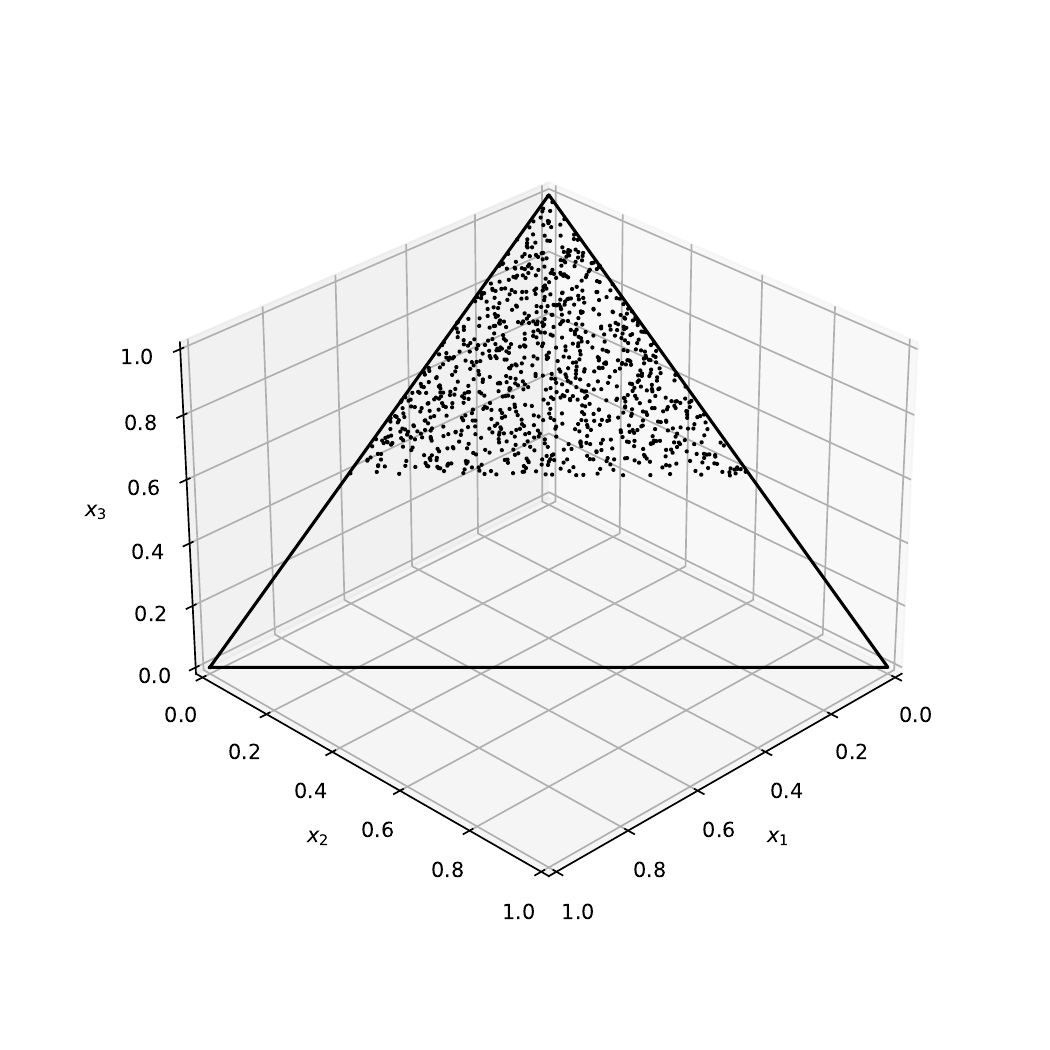}
        \caption[]%
        {{\small Inequality constraint: $x_1+x_2\geq 0.6$}}    
        \label{fig:drsc_inequality}
    \end{subfigure}
    \hfill
    \begin{subfigure}[b]{0.44\textwidth}   
        \centering 
        \includegraphics[width=\textwidth]{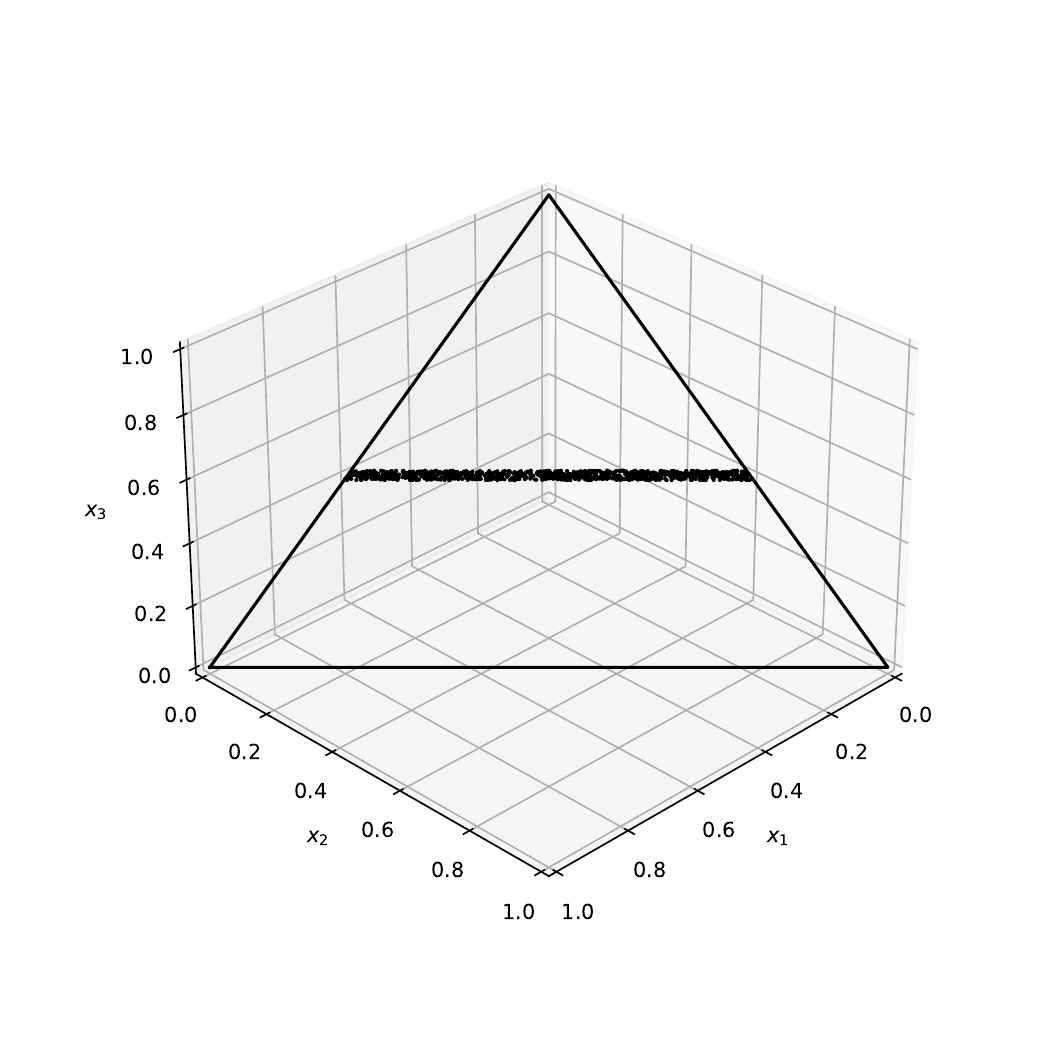}
        \caption[]%
        {{\small Equality constraint: $x_1+x_2= 0.6$}}    
        \label{fig:drsc_equality}
    \end{subfigure}
    \vskip\baselineskip
    \begin{subfigure}[b]{0.44\textwidth}   
        \centering 
        \includegraphics[width=\textwidth]{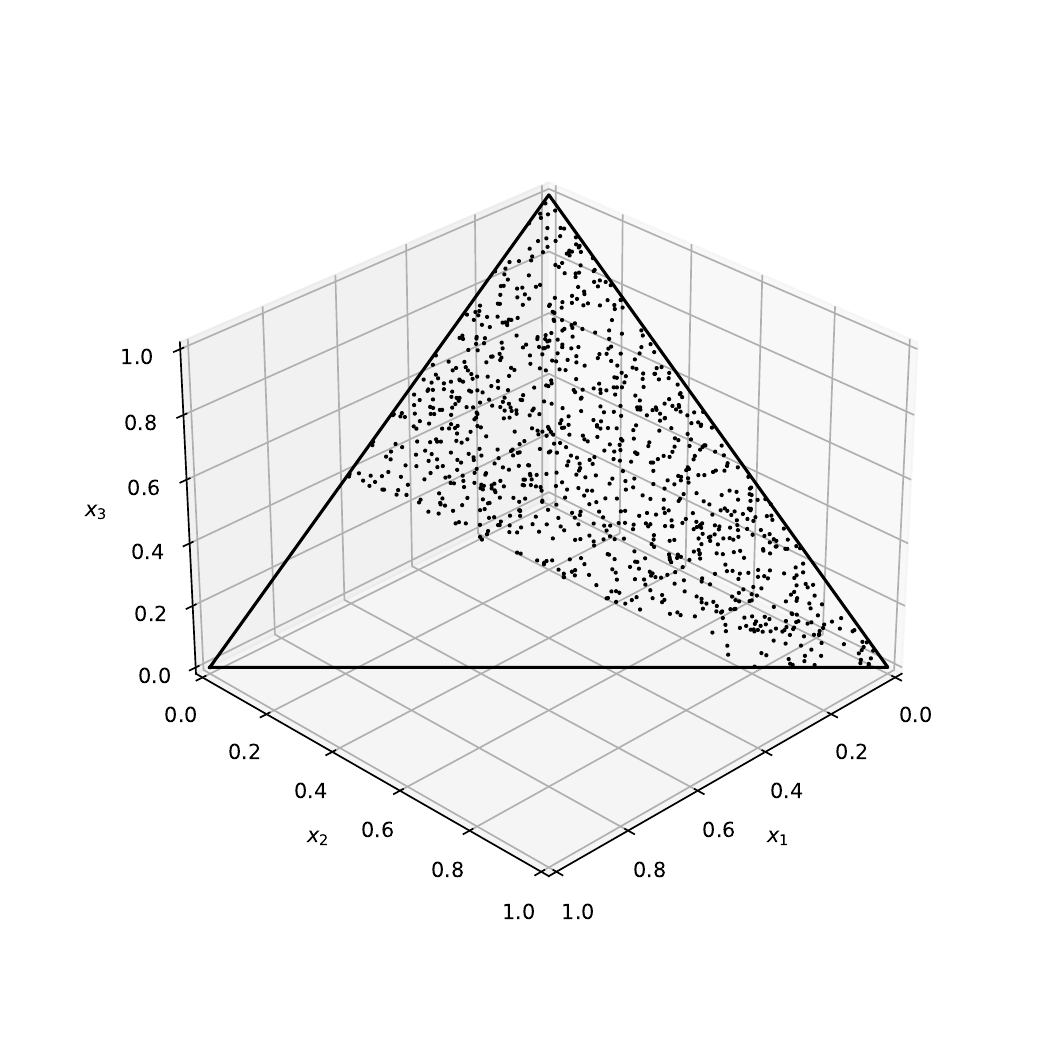}
        \caption[]%
        {{\small Inequality constraint: $x_1+0.5 x_2 \leq 0.6$}}    
        \label{fig:drsc_general_inequality}
    \end{subfigure}
    \hfill
    \begin{subfigure}[b]{0.44\textwidth}   
        \centering 
        \includegraphics[width=\textwidth]{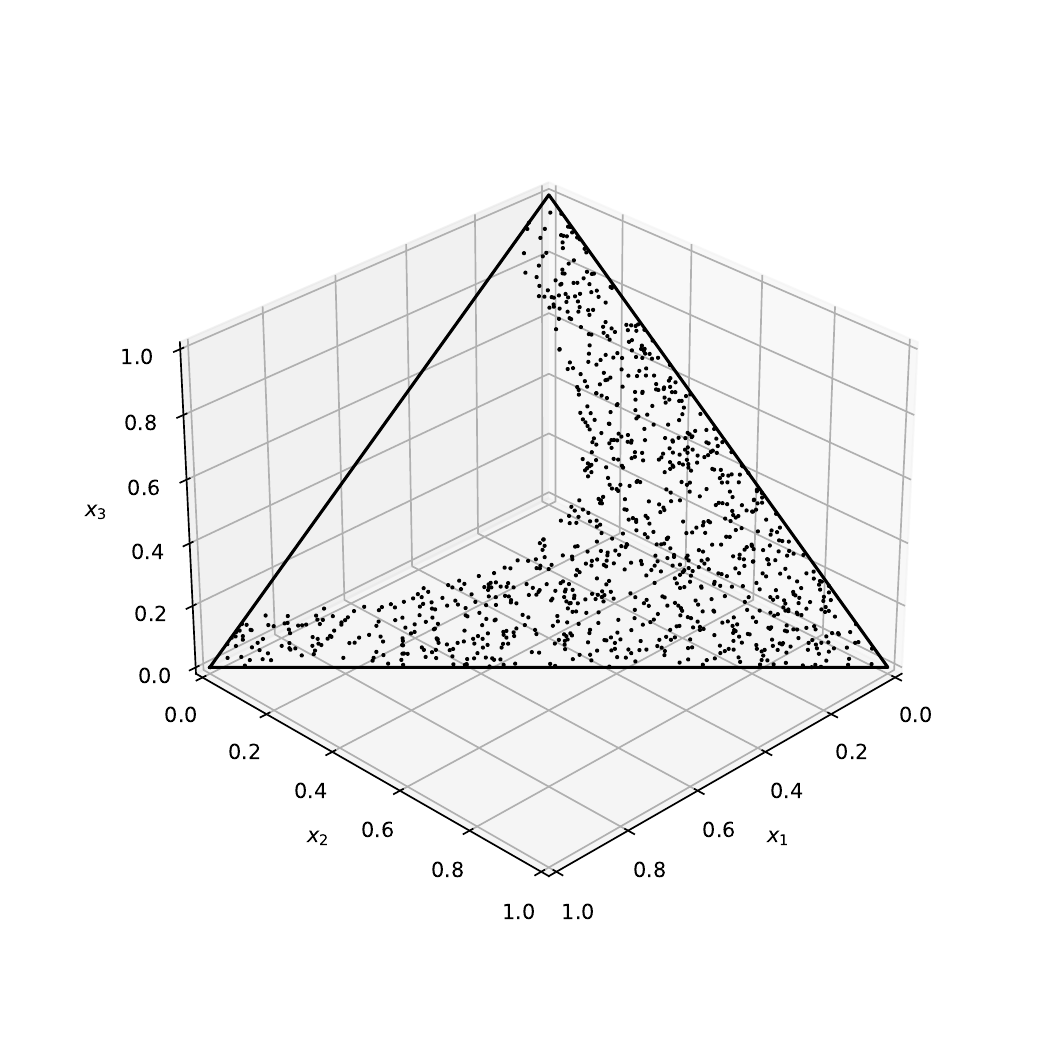}
        \caption[]%
        {{\small Nonlinear inequality: $x_1 x_2-0.1\leq 0$.}}    
        \label{fig:standardsimplex_withquadraticinequalityconstraints}
    \end{subfigure}
    \caption[]
    {\small Samples from standard simplices with constraints. The samples are obtained using 1000 draws from the DRSC algorithm.} 
    \label{fig:drsc}
\end{figure*}

For the DRSC algorithm we set the feasibility parameters to $\epsilon_{ineq}=10^{-3}$ and $\epsilon_{eq}=10^{-2}$. We show that the DRSC algorithm can generate a uniform distribution of vectors under different conditions. In \Cref{fig:drsc}, we illustrate a generated sample of 1000 $3$-dimensional vectors in the standard simplex. In each subfigure, we consider an additional constraint. In \Cref{fig:drsc_inequality} the inequality $x_1+x_2\geq 0.6$ is used. \Cref{fig:drsc_equality} models an equality constraint $x_1+x_2= 0.6$ and \Cref{fig:drsc_general_inequality} showcases a non-symmetric constraint $x_1+0.5 x_2 \leq 0.6$. In \Cref{fig:standardsimplex_withquadraticinequalityconstraints}, we generate vectors satisfying a nonlinear constraint $x_1 x_2-0.1\leq 0$. We observe that when adding linear and nonlinear constraints, the vectors in the figures indeed seem to be drawn from a uniform distribution.

These experiments demonstrate the capabilities of the DRSC algorithm to handle various types of constraints, but do not reflect typical use cases. For instance, in a $3$-dimensional standard simplex the constraint $x_1+x_2\geq 0.6$ is equivalent to $x_3\leq 0.4$. However, in higher dimensions such a substitution is no longer possible.
\FloatBarrier

\section{Conclusion and discussion}\label{section7:conclusion}
We consider the problem of generating a uniform distribution of vectors with fixed sum satisfying both linear and nonlinear constraints. \cite{griffin2020} proposed the Dirichlet-Rescale (DRS) algorithm to generate vectors with individual lower and upper bounds. However, we demonstrate that the vectors obtained using the DRS algorithm are not necessarily generated from a uniform distribution. This happens for instance when multiple bounds may be violated at the same time. 

In practice we may expect that the vectors generated by the DRS algorithm from \cite{griffin2020} are close to a uniform distribution. This observation is supported by statistical test in which we need a large number of vectors to detect that the sample is not uniformly distributed. 

We introduce the Dirichlet-Rescale-Constraints (DRSC) algorithm, which handles both linear and nonlinear constraints, while ensuring that each vector within the valid region is drawn with an equal probability. The main idea of the DRSC algorithm is to find large regular induced simplices, which do not overlap with the other induced simplices.

A limitation of the current DRSC algorithm is the number of restarts  when a vector violates a (non)linear constraint but does not lie within a regular induced simplex. To minimise the number of restarts one should cover the infeasible region by identifying large, possibly irregular simplices. However, identifying large (irregular) induced simplices presents a challenging optimisation problem.

\renewcommand{\bibname}{References} 
\phantomsection
\addcontentsline{toc}{section}{References}
\bibliographystyle{apalike}
\bibliography{sample} 

\newpage
\phantomsection
\addcontentsline{toc}{section}{Appendix}



\setcounter{table}{0}
\renewcommand{\thetable}{A\arabic{table}}
\renewcommand{\thesubsection}{\Alph{subsection}} 

\end{document}